\documentclass[12pt,leqno]{amsart}
   
\usepackage[top=2.5 cm,bottom=2 cm,left=2.5 cm,right=2 cm]{geometry}
\usepackage[utf8]{inputenc}
\usepackage[dvipsnames]{xcolor}
\usepackage{tikz}
\usepackage{pgfplots}
\usepackage{footmisc}
\pgfplotsset{compat=newest}
\usepgfplotslibrary{fillbetween}
\usetikzlibrary{patterns}
\usepackage{bbm}
\usepackage{esint}
\usepackage{graphicx}
\usepackage{hyperref}
\usepackage{color}
\usepackage{amsfonts}
\usepackage{psfrag}
\newcounter{stepnb}

\newtheorem{theorem}{Theorem}[section]

\newtheorem{definition}{Definition}[section]

\numberwithin{equation}{section}

\newcommand{\R}{\mathbb{R}}

\newcommand{\ee}{\varepsilon}
\newcommand{\mf}{\mathbf}

\usepackage{pgfplots}

\newcommand{\be}{\begin{equation}}
\newcommand{\eq}{\end{equation}}

\begin{document}
\dedicatory{Dedicated to Professor  Gui-Qiang Chen on the occasion of his 60th birthday}
\title[]{An informal account of recent results on initial-boundary value problems for systems of conservation laws}

\author[L.~V.~Spinolo]{Laura V.~Spinolo}
\address{L.V.S. CNR-IMATI, via Ferrata 5, I-27100 Pavia, Italy.}
\email{spinolo@imati.cnr.it}
\author[F.~Ancona]{Fabio Ancona}
\address{F.A. Dipartimento di Matematica ``Tullio Levi-Civita", Universit\`a degli Studi di Padova, Via Trieste 63,
35131 Padua, Italy}
\email{ancona@math.unipd.it}
\author[A.~Marson]{Andrea Marson}
\address{A.M. Dipartimento di Matematica ``Tullio Levi-Civita", Universit\`a degli Studi di Padova, Via Trieste 63,
35131 Padua, Italy}
\email{marson@math.unipd.it}
\maketitle
{
\rightskip .85 cm
\leftskip .85 cm
\parindent 0 pt
\begin{footnotesize}

{\sc Abstract.} This note aims at providing a rather informal and hopefully accessible overview of the fairly long and technical work~\cite{AFL}. In that paper, the authors established new global-in-time existence results for admissible solutions of nonlinear systems of conservation laws defined in domains with boundaries. The main novelty in~\cite{AFL} is that the solution is constructed by taking into account the underlying viscous mechanism, which is relevant because, in the case of initial-boundary value problems, different viscous approximations yield in general different limits. This note will frame the analysis of~\cite{AFL} in the relevant context, compare the main result with the previous existing literature, and touch upon the most innnovative technical points of the proof. 




\end{footnotesize}
}

\section{Introduction} 
The goal of this note is to provide an handwaving overview of some recent progress on the analysis of the initial-boundary value problem for nonlinear systems of conservation laws in the form 
\begin{equation} 
\label{e:claw}
      \mf  g (\mf v)_t +  \mf f(\mf v)_x = \mf 0.
\end{equation}
In the above expression, the unknown $\mf v$ depends on the variables $(t, x)$, with $x$ one-dimensional, and attains values in $\R^N$. The functions $\mf g, \, \mf f: \R^N \to \R^N$ are smooth and satisfy suitable assumptions that we touch upon in the following. The exposition in this note is mainly based on the recent paper~\cite{AFL}, where~\eqref{e:claw} is coupled with the viscous approximation 
\begin{equation} 
\label{e:vclaw}
       \mf g (\mf v^\ee)_t + \mf f(\mf v^\ee)_x = \ee \Big( \mf D(\mf v^\ee) \mf v^\ee_x \Big)_x. 
       \end{equation}
In the above expression, $\mf D$ is a positive semi-definite $N \times N$ matrix depending on the physical model under consideration. We discuss in the following the precise assumptions we impose on $\mf D$, for the time being we only mention that we rely on the analysis in the fundamental works by Kawashima and Shizuta, see~\cite{KawashimaShizuta1}.  In particular, our assumptions apply in the case of the compressible Navier-Stokes (or Navier-Stokes-Fourier) equations and the viscous MHD equations, which are the most natural choices for approximating the Euler and inviscid MHD equations, respectively. From the analytical standpoint, a very relevant feature of these cases (and of most, if not all, the physically relevant cases) is that the matrix $\mf D$ is singular, so that~\eqref{e:vclaw} is a \emph{mixed hyperbolic-parabolic} system. This accounts for severe technical challenges that we touch upon in the following. 

Coupling~\eqref{e:claw} with the underlying viscous mechanism is especially important in the case of initial-boundary value problems because, as we will discuss in the following, different choices of $\mf D$ in~\eqref{e:vclaw} yield in general different solutions of~\eqref{e:claw} in the vanishing $\ee$ limit. The main result of~\cite{AFL}, which is also Theorem~\ref{t:main} below, establishes global-in-time existence of admissible solutions of~\eqref{e:claw} consistent with the underlying viscous mechanism~\eqref{e:vclaw}.

The outline of this note is as follows. In \S\ref{s:bi} we discuss  the main motivation for the analysis in~\cite{AFL} by providing some background information and focusing in particular on the viscous approximation of initial-boundary value problems. In \S\ref{s:mr} we state the main result of~\cite{AFL}, extensively comment on it, and provide a very high level overview of the proof. In \S\ref{s:overview} we compare the main result of~\cite{AFL} with the existing literature, and finally in \S\ref{s:technical} we explain what we feel are, from the technical standpoint, the most interesting points of the proof of Theorem~\ref{t:main}.

\section{Background information} \label{s:bi}
\subsection{Nonlinear systems of conservation laws in one space variable, and their viscous approximation}
The archetype of~\eqref{e:claw} are the celebrated compressible Euler equations modeling the dynamics of an ideal and compressible fluid. In this case, $N=3$ and the component of the unknown $\mf v$ represent the fluid density, velocity, and internal energy. Other famous examples include the inviscid Magneto-Hydro-Dynamic (MHD) equations describing the propagation of plane waves in an electrically charged, compressible and ideal fluid. 

Note however that the present note, despite using the Euler and MHD equation as a guiding thread for the exposition, aims at discussing results that apply to \emph{general} systems of conservation laws. 
In this framework and in the case of the Cauchy problem, Glimm~\cite{Glimm} established existence of global-in-time, admissible solutions of~\eqref{e:claw} by relying on a suitable random choice approximation method now named after him. 
Uniqueness results were obtained in a series of works by Bressan \emph{et al.}, see~\cite{Bressan} and also the very recent contributions~\cite{BressanDeLellis,BressanGuerra}. A fundamental tool in the proof of the above uniqueness results is the analysis of the so-called \emph{wave front-tracking scheme}, a suitable approximation algorithm through piecewise constant functions that we briefly discuss in \S\ref{ss:wft} below. Note furthermore, and more importantly, that this well-posedness theory requires that the total variation of the initial data is sufficiently small, with the smallness threshold depending on the specific system under consideration. Well-posedness results for data of large total variation are only available for very specific systems, like those in the so-called \emph{Temple class}~\cite{Temple}, or in special perturbative regimes, see~\cite{Lewicka}. This is a somehow unavoidable drawback of the theory since there are explicit examples~\cite{Jen} showing the if the total variation of the initial datum is finite, but large, then the admissible solution can blow-up in finite time. The blow-up system in~\cite{Jen} is not physical (i.e., it does not admit strictly convex entropies), but the recent contributions~\cite{BCZ2,BCZ} provide strong indications of a possible finite time blow-up of the total variation in the case of the so-called $p$-system of isentropic gas dynamics in Lagrangian coordinates, arguably one of the simplest systems with a clear physical meaning.  

As mentioned before, in this note we will recover the inviscid system~\eqref{e:claw} as the limit of the underlying viscous mechanism~\eqref{e:vclaw}. Indeed, the viscous system~\eqref{e:vclaw} formally boils down to~\eqref{e:claw} as the parameter $\ee$ vanishes. However, despite considerable efforts and advances,  
the analytical understanding of the vanishing viscosity limit is still badly incomplete and  a general convergence proof is missing. In the \emph{artificial viscosity} case, that is when the matrix $\mf D$ is the identity, the right-hand side of~\eqref{e:vclaw} boils down to $\ee  \mf v^\ee_{xx}$ and a very complete result was established by Bianchini and Bressan in~\cite{BianchiniBressan}. In the \emph{physical case} the situation is much less understood and fewer results are available. Among them, the remarkable work by G.Q. Chen and Perepelitsa~\cite{ChenPerepelitsa} applies to system of compressible and barotropic gas dynamics, that is to the $2\times2$ Navier-Stokes equations in Eulerian coordinates, and relies on the compensated compactness techniques discussed in the pivotal works by DiPerna~\cite{DiPerna1,DiPerna2}. We also refer to the very recent contribution by G. Chen, Kang and Vasseur~\cite{ChenKangVasseur}, where the authors manage to recover small $BV$ solutions of the system of compressible and barotropic inviscid gas dynamics as strong limits of solutions of the physical viscous approximation, that is the $2\times2$ Navier-Stokes equations. 

\subsection{On the viscous approximation of nonlinear systems of conservation laws in domains with boundaries} \label{ss:va}
In the case of initial-boundary value problem, the analysis of the vanishing viscosity limit of~\eqref{e:vclaw} is further complicated by the presence of boundary layer phenomena, that account for a loss of boundary condition and for a sharp transient behavior of the solution of the viscous system close to the domain boundary. Nevertheless, partial convergence results have been obtained in some works, see for instance~\cite{AnconaBianchini,ChenFrid,Gisclon,GrenierRousset,JosephLeFloch,Rousset,Spinolo}. See also~\cite{Serre1,Serre2} for a general introduction to initial-boundary value problems for conservation laws. 

The signature feature of the viscous approximation of the initial-boundary value problems, and a striking difference with the Cauchy problem case, is that the vanishing $\ee$ limit of~\eqref{e:vclaw} in general \emph{depends} on $\mf D$ and changes as $\mf D$ changes, see~\cite{Gisclon}. Remarkably, this happens even in the most elementary linear case: consider the linear viscous system
\be \label{e:linear}
   \mf v^\ee_t + \mf A \mf v^\ee_x = \ee \mf D \mf v^\ee_{xx},
\eq
where $\mf A$ and $\mf D$ are (constant) $N\times N$ matrices. To highlight the heart of the matter and avoid some technicalities, let us assume that $\mf A$ is symmetric and invertible and that $\mf D$ is symmetric and positive definite. We assume that~\eqref{e:linear} is defined on the domain $x > 0$, $t \ge 0$ and we  couple~\eqref{e:linear} with the Riemann-type data
\be
\label{e:rie}
       \mf v^\ee (0, x) = \mf v_0, \qquad \mf v^\ee (t, 0) = \mf v_b,
\eq
where $\mf v_0 $ and $\mf v_b$ are given states in $\R^N$. In this case, fairly standard energy-type estimates ensure that $\mf v^\ee$ converges weakly in $L^2_{\mathrm{loc}} (\R_+ \times \R)$ to a weak solution $\mf v$ of the linear advection equation 
\be
\label{e:linearcl}
        \mf v_t + \mf A \mf v_x =0
\eq
which has the following structure: it is piecewise constant, with contact discontinuities occurring along $p$ half-lines departing from the origin. Here $p$ denotes the number of \emph{positive} eigenvalues of $\mf A$. The solution $\mf v$ attains the boundary datum $\mf{\bar v}: = \mf v(t, 0)$, which does not depend on $t$ since $\mf v$ is self-similar. Most importantly, the relation between $\mf{\bar v}$ and the boundary datum $\mf v_b$ imposed on the viscous approximation is the following: there is a so-called \emph{boundary layer} $\mf w: [0, + \infty[ \to \R^N$ such that 
\be \label{e:bllin}
     \left\{
      \begin{array}{ll}
             \mf A \mf w' = \mf D \mf w'' \\
             \mf w(0) = \mf v_b, \quad \lim_{y \to + \infty} \mf w(y) = \mf{\bar v}. 
      \end{array}
     \right.
\eq
Note that the equation at the first line of the above system dictates that the boundary layer $\mf w$ is, up to a change of variables, a steady solution of the viscous equation~\eqref{e:linear}. Fairly classical results on linear systems imply that the boundary value problem~\eqref{e:bllin} admits a solution if and only if $\mf v_b - \mf{\bar v}$ belongs to the so-called \emph{stable subspace}, that is to the generalized eigenspace of $\mf D^{-1} \mf A$ associated to eigenvalues with negative real part. Since this eigenspace obviously \emph{depends} on $\mf D$, so does the boundary value $\mf{\bar v}$ and hence the limit solution $\mf v$. 

The above elementary argument describes the very basic mechanism yielding the dependence of the vanishing viscosity approximation on the underlying viscous system. Note that it is a mechanism, basically due to boundary layer phenomena, that acts in the initial-boundary case only, and indeed in the case of the Cauchy problem it can be shown that, at least in suitable small total variation regimes, the limit of the vanishing viscosity~\eqref{e:vclaw}, if any, does not depend on $\mf D$, see~\cite{Bianchini}. As a side remark, the dependence on the underlying viscous mechanism of solutions of conservation laws defined on domains with boundary has also very relevant consequences from the numerical viewpoint. Very loosely speaking, this is due to the fact that most numerical schemes for conservation laws contain the so-called \emph{numerical viscosity}. In standard numerical schemes, the numerical viscosity is modeled upon the \emph{artificial viscosity}, namely the viscosity mechanism obtained by choosing as viscosity matrix $\mf D$ the identity. As a consequence, standard numerical scheme provide an approximation of the inviscid limit of the artificial viscosity, which, in the case of initial-boundary value problems, differs from the physically relevant solution, which is the limit of the 
\emph{physical} viscosity. To fix this issue, one can introduce new numerical schemes where the numerical viscosity is modeled upon the physical viscosity, see~\cite{MishraSpinolo} for related numerical experiments. 
\section{Main result} \label{s:mr}
In this section we state the main result of~\cite{AFL}, which is Theorem~\ref{t:main} below. More precisely, in \S\ref{ss:bc} we provide the formulation of the boundary condition assigned on~\eqref{e:claw}, a non trivial issue in view of the discussion in \S\ref{ss:va}. In \S\ref{ss:ge} we provide the precise statement of Theorem~\ref{t:main}, and make some comments. Finally, in \S\ref{ss:over} we discuss the proof of Theorem~\ref{t:main} at a fairly high level, and refer to \S\ref{s:technical} for more technical comments. 
\subsection{Formulation of the initial-boundary value problem} \label{ss:bc}
Going back to the analytic study of the initial-boundary value problem for~\eqref{e:claw}, from the discussion in \S\ref{ss:va} and in particular from the considerations on the linear case~\eqref{e:linearcl} we infer that  the key point to understand the dependence of the viscous limit on the viscosity matrix $\mf D$ is to unveil the relation between the boundary datum imposed on the \emph{viscous} system~\eqref{e:vclaw} and the trace of the solution of the inviscid conservation law~\eqref{e:claw}. 
Towards this end, we first of all point out that, owing to the singularity of the matrix $\mf D$, the initial-boundary value problem for~\eqref{e:vclaw} is in general overdetermined if one imposes a \emph{full} boundary condition like $\mf v^\ee (t, 0) = \mf v_b (t)$.  To see this, let us consider a specific example: the first line of both the Navier-Stokes and the viscous MHD equation written in Eulerian coordinates expresses mass conservation and reads 
$$
    \partial_t \rho + \partial_x [\rho u]=0,
$$
where $\rho$ and $u$ denote the fluid density and velocity, respectively. If $u<0$ at the boundary then assigning both the boundary and the initial condition on $\rho$ would yield an overdetermined problem. Going back to the general case, under fairly resonable assumptions on the mixed hyperbolic-parabolic system one can introduce a slightly involved formulation of the boundary condition to restore well-posedness, namely one couples~\eqref{e:vclaw} with the initial and boundary conditions 
\be 
\label{spinolo_e:datavclaw}
   \mf v^\ee(0, \cdot) = \mf v_0, \qquad \boldsymbol{\widetilde \beta} (\mf v^\ee (\cdot, 0), \mf v_b) = \mf 0_{N}. 
\eq
The precise definition of the function $ \boldsymbol{\widetilde \beta}: \R^N \times \R^N \to \R^N$ is provided in \cite[\S2.2]{AFL}, but the very basic idea underpinning the construction of $\boldsymbol{\widetilde \beta}$ is imposing a full boundary condition on the \emph{parabolic} component of $\mf v^\ee$, and a boundary condition along the characteristic fields of the \emph{hyperbolic} component entering the domain. In particular, going back to specific example of the Navier-Stokes equations, prescribing $\boldsymbol{\widetilde \beta} (\mf v^\ee (\cdot, 0), \mf v_b) = \mf 0_{N}$ means that we \emph{always} assign the values of the fluid velocity $u$ and of the internal energy at the boundary, and that if $u>0$ at the domain boundary we also assign the value of the fluid density, otherwise we do not.  

With the above notation in place, we can now discuss the boundary condition we impose on the inviscid system~\eqref{e:claw}. We proceed step by step and at first we introduce an additional assumption that considerably simplifies the matter. We term the domain boundary $x=0$ \emph{not characteristic} if all the eigenvalues of the jacobian matrix $\mf D \mf f (\mf u)$ are bounded away from $0$, for every $\mf u \in \R^N$. In other words, in the non characteristic boundary case we can easily discriminate between the characteristic field of~\eqref{e:claw} \emph{entering} and \emph{leaving} the domain at the domain boundary. In this framework, we impose the initial and boundary condition on~\eqref{e:claw} by requiring that 
\be
\label{e:ibvp}
         \mf v(0, \cdot) = \mf v_0 \; \; \text{a.e. on $\R_+$}, \qquad 
         \mf v(\cdot, 0) \sim_{\mf D} \mf v_b  \; \; \text{a.e. on $\R_+$},
\eq
where the relation $\sim_{\mf D}$ is defined as follows: given $\mf{\bar v}, \mf v_b \in \R^N$ we say that  $\mf{\bar v}  \sim_{\mf D} \mf v_b$ if 
there is a boundary layer $\mf w: \R_+ \to \R^N$ such that
\be
\label{e:bl2}
\left\{
\begin{array}{ll}
          \mf D(\mf w) \mf w' = \mf f(\mf w) - \mf f(\mf{\bar v}) \\
           \boldsymbol{\widetilde \beta} (\mf w(0), \mf v_b) =\mf 0_N, \quad \lim_{y \to + \infty} \mf w(y)=\mf{\bar v}. 
\end{array}
\right.
\eq
Keeping in mind that  $\boldsymbol{\widetilde \beta} (\mf w(0), \mf v_b) =\mf 0_N$ is the way we assign the boundary condition on the viscous system~\eqref{e:vclaw}, we conclude that~\eqref{e:bl2} is, up to space integration, the exact nonlinear analogous of~\eqref{e:bllin}. 

We now move towards the most general case and take into account the possibility that an eigenvalue of the jacobian matrix $\mf D \mf f \mf (\mf u)$ attains the value $0$: this is usually refered to as the \emph{boundary characteristic case}. Note that for instance the boundary $x=0$ is characteristic for the compressible Euler equations written in Eulerian coordinates if the fluid velocity either vanishes, or is close (in modulus) to the sound speed. The boundary $x=0$ is always characteristic for the compressible Euler equations written in Lagrangian coordinates, as the second eigenvalue vanishes identically.

Compared to the treatment of the non-characteristic case, the analysis of the boundary characteristic case involves additional severe technical challenges that we touch upon in the following. For the time being, we mention that in the boundary characteristic case we have to take into account the possibility of a shock or contact discontinuity that sits exactly at the domain boundary. With the above preliminary considerations in place, we can now introduce the general definition of the relation~$\sim_{\mf D}$. For simplicity, from now on we focus on the case where the characteristic field associated to the vanishing eigenvalue is either genuinely nonlinear or linearly degenerate: in the following, we term this field \emph{boundary characteristic field}. Note that in many physically relevant systems, like the compressible Euler equations written in both Eulerian and Lagrangian coordinates, every vector field (and henceforth, also the boundary characteristic field, if any) is either genuinely nonlinear or linearly degenerate. 
\begin{definition}
\label{d:equiv} Given system~\eqref{e:vclaw} and $\mf {\bar v}$, $\mf{ v}_b \in \R^N$, we say that ``$\mf{\bar v} \sim_{\mf D} \mf{ v}_b$'' if there is 
$\mf{\underline v} \in \R^N$ such that the following conditions are both satisfied:
\begin{itemize}
\item[i)] $\mf f(\mf{\bar v}) = \mf f(\mf{\underline v})$; if the boundary characteristic field is genuinely nonlinear, we also require that the $0$-speed shock between  $\mf{\bar v}$ (on the right) and $\mf{\underline v}$ (on the left) is Lax admissible\footnote{if the boundary characteristic field is linearly degenerate, then the discontinuity is a contact discontinuity and the Lax admissibility condition is automatically satisfied}; 
\item[ii)] there is a so-called ``boundary layer'' $\mf w: \R_+ \to \R^N$ such that
\be
\label{e:bl}
\left\{
\begin{array}{ll}
          \mf D(\mf w) \mf w' = \mf f(\mf w) - \mf f(\mf{\underline v}) \\
           \boldsymbol{\widetilde \beta} (\mf w(0), \mf v_b) =\mf 0, \quad \lim_{y \to + \infty} \mf w(y)=\mf{\underline v}. 
\end{array}
\right.
\eq
\end{itemize}
\end{definition}
Note that, if we apply Definition~\ref{d:equiv} in the non characteristic boundary case, if $\mf{\bar v}$ is confined in a small enough neighborhood of $\mf{\underline v}$ the identity $\mf f(\mf{\bar v}) = \mf f(\mf{\underline v})$ implies by the Local Invertibility Theorem that $\mf{\bar v}= \mf{\underline v}$, and hence~\eqref{e:bl} boils down to~\eqref{e:bl2}. 
 \subsection{Global-in-time existence of admissible solutions} \label{ss:ge}
With Definition~\ref{d:equiv} in place, we can now state the main result of~\cite{AFL}. 
\begin{theorem}
\label{t:main} Assume that $\mf g$, $\mf f$, $\mf D$ satisfy Hypotheses $1, \cdots, 5$ in~\cite[\S2]{AFL}, and fix $\mf v^\ast \in \R^N$; then there is a constant $\delta^\ast>0$ only depending on the  functions $\mf g$, $\mf f$, $\mf D$ in system~\eqref{e:vclaw} such that the following holds. If $\mf v_0, \mf v_b \in BV(\R_+)$ satisfy
\be \label{e:hp}
     \mathrm{TotVar} \ \mf v_0 + \mathrm{TotVar} \ \mf v_b + |\mf v_0 (0^+) - \mf v_b (0^+)| \leq \delta^\ast, \quad |\mf v_0 (0^+) - \mf v^\ast| \leq \delta^\ast
\eq
then there is a global-in-time, Lax admissible distributional solution $\mf v \in BV_{\mathrm{loc}} (\R_+ \times \R_+)$ of~\eqref{e:claw} satisfying the initial and boundary conditions~\eqref{e:ibvp} in the sense of Definition~\ref{d:equiv}. Also, if system~\eqref{e:claw} admits a convex entropy then the solution we construct is  entropy admissible.   
\end{theorem}
In the statement of the above result, $\mathrm{TotVar}\ \mf v$ and $BV$ denote the total variation of the function $\mf v$ and the space of bounded total variation functions, respectively. Also, $\mf v_0 (0^+)$ and $\mf v_b (0^+)$ denotes the right limit of $\mf v_0$ and $\mf v_b$ at $0$, respectively. Note furthermore that, as a function of total variation, the function $\mf v$ has a well-defined trace at $x=0$. 

We refer to \S\ref{s:overview} below for a more detailed comparison between Theorem~\ref{t:main} and previous related results. Here we just point out that the main novelty of Theorem~\ref{t:main} is that, in the boundary characteristic case, it provides what is to the best of our knowledge the first global existence result for solutions of nonlinear systems of conservation laws~\eqref{e:claw} consistent with the underlying viscous mechanism~\eqref{e:vclaw}. This consistency is encoded in the boundary condition of Definition~\ref{d:equiv}. As we discussed before, enclosing information on the underlying viscous mechanism in the definition of solution is fundamental, because different viscous approximation yield different inviscid limits in the case of initial-boundary value problems.  

Rather than providing the exact statement of the hypotheses of Theorem~\ref{t:main} concerning $\mf g$, $\mf f$, $\mf D$, which would require the introduction of some heavy notation, we now informally discuss them. First of all, we stress that Hypotheses $1, \cdots, 5$ in~\cite[\S2]{AFL} are all satisfied by our archetipycal physical examples, that is by the compressible Euler and MHD equations, written in both Eulerian and Lagrangian coordinates. More in detail, Hypothesis 1 states that  system~\eqref{e:vclaw} is, up to a change of the dependent variables, in the \emph{normal form}, in the Kawashima-Shizuta sense~\cite{KawashimaShizuta1}, that is it allows for a specific block decomposition. Hypothesis 2 states that system~\eqref{e:vclaw} satisfies the so-called \emph{Kawashima-Shizuta condition}, which heuristically speaking is a coupling condition that rules out the possibility of singling out in~\eqref{e:vclaw} a purely hyperbolic component. Hypothesis 3 is a fairly standard strict hyperbolicity assumption. Hypothesis 4 is also fairly common in the conservation laws analysis framework and dictates that every characteristic field of~\eqref{e:claw} is either linearly degenerate or genuinely nonlinear. Finally, Hypothesis 5 requires that system~\eqref{e:vclaw} satisfies some slightly technical conditions introduced in the previous works~\cite{BianchiniSpinoloARMA,BianchiniSpinolo} whose meaning is, loosely speaking, to ensure  that the boundary layers system~\eqref{e:bl} can be written in a tractable form, a nontrivial requirement given the singularity of the viscosity matrix $\mf D$. 

Concerning the small total variation assumption~\eqref{e:hp}, this is obviously highly restrictive, but it is basically a necessary price to pay if the goal is to deal with fairly general systems, since as mentioned before there are by now very strong indications that, if the smallness assumptions fails, then finite time total variation blow-up may occur even in the case of relativey simple examples
like the $p$-system, see~\cite{BCZ}.
 
As a last comment, we point out that Theorem~\ref{t:main} is a global existence result only: we are fairly confident one could establish uniqueness results by relying on the Standard Riemann Semigroup approach \emph{\`a la} Bressan, see also the related work~\cite{AmadoriColombo}. However, we decided to postpone the uniqueness proof to the future given that~\cite{AFL} is already fairly long and technical.

\subsection{A very non technical overview of the proof of Theorem~\ref{t:main}} \label{ss:over}
We now briefly and informally comment on the proof of Theorem~\ref{t:main}, and we  
refer to \S\ref{s:technical} below for some more detailed technical remarks concerning the main novelties of our construction. In a nutshell, the proof of Theorem~\ref{t:main} relies on the introduction of a new wave front-tracking algorithm, and the most delicate points of our analysis stem from the fact that we deal with the boundary characteristic case. At a very high level, the presence of a possibly vanishing boundary characteristic field implies that one cannot make a clear distinction between waves that are entering and leaving the domain $x>0$ at the domain boundary. In particular, one should in principle take into account the possibility of waves of the boundary characteristic family that are bounced back and forth at the boundary, a behavior that could lead to the finite time blow-up of the total variation and to the breakdown of the wave front-tracking algorithm. At a more technical level, our analysis relies on a detailed description of the structure of the boundary layers $\mf w$ solving~\eqref{e:bl} that was provided in the previous works~\cite{BianchiniSpinoloARMA,BianchiniSpinolo}. If the boundary is characteristic, the boundary layers may have a component lying on a suitable center manifold: by slightly perturbing this component one can obtain waves (rarefaction waves, shocks or contact discontinuities) of the boundary characteristic family that travel with very small albeit positive speed 
and hence enter the domain.  This possibility severely complicates our analysis, and requires the introduction of several new ideas to obtain useful estimates on the total variation increase at times where a boundary characteristic wave hits the boundary, estimates that in the wave front-tracking jargon are called \emph{interaction estimates}. We come back to this point in \S\ref{s:technical}. 

Note that relying on wave front-tracking techniques often accounts for fairly long and technical proofs, but also carries several advantages and paves the way for further developments. In particular, the introduction of a suitable wave front-tracking algorithm is pivotal to the proof of the uniqueness results for~\eqref{e:claw} due to Bressan and collaborators, see~\cite{Bressan}. Also, it is known that for instance a $BV$ solution $\mf v$ of the conservation law~\eqref{e:claw} recovered as the limit of a wave front-tracking algorithm enjoys better regularity properties than a generic function of bounded variation, see for instance the analysis in~\cite[Ch. 10]{Bressan}.

\section{Comparison with previous results}\label{s:overview}
\subsection{The scalar case}
To start our little overview of global-in-time existence results for solutions of~\eqref{e:claw} defined on domains with boundary we briefly touch up the scalar case $N=1$. In this case, and for conservation laws in several space dimensions, the paper~\cite{BLRN} by Bardos, Le Roux and Nedelec extends the analysis of the milestone work of Kru{\v{z}}kov
to the case of domains with boundaries. In particular, it establishes convergence of the vanishing viscosity approximation and uniqueness of (a suitable notion of) entropy admissible solution. 
\subsection{Global in time existence results \emph{via} Glimm-type schemes}
Moving towards systems, i.e towards the case $N>1$, one of the very first works dealing with initial-boundary value problems is the paper by Nishida and Smoller~\cite{NishidaSmoller}. It establishes global-in-time existence results for the so-called \emph{piston problem}, namely the initial-boundary value problem for the $p$-system of isentropic gas dynamics in Lagrangian coordinates with the values of the fluid velocity prescribed at the domain boundary. The analysis in~\cite{NishidaSmoller} relies on the introduction of a suitable modification of the approximation algorithm introduced in~\cite{Glimm} and by-now known as \emph{the Glimm scheme}.  Also based on the analysis of a Glimm-type algorithm is the paper by Liu~\cite{LiuIBVP} concerning the full compressible Euler equations written in Lagrangian coordinates. This is a system of $3$ equations where the second eigenvalue of the jacobian matrix of the flux is identically $0$, which implies that the boundary  is characteristic. Note, however, that the fact that the characteristic eigenvalue vanishes identically (and not only at some point) considerably simplifies the treatment of the boundary characteristic case as it still allows for a clear discrimination between wave entering and leaving the domain at the boundary, with the waves of the second family sitting exactly at the domain boundary and tangent to it. In~\cite{LiuIBVP}, the author assigns the value of either the fluid velocity or the pressure at the domain boundary.

The most general contributions that  rely on the introduction of Glimm-type schemes are the PhD thesis of Goodman~\cite{Goodman} and the work by  Sabl{\'e}-Tougeron~\cite{SableT}. In particular, Goodman establishes global-in-time existence of admissible solutions for general systems in the non-characteristic boundary case. Sabl{\'e}-Tougeron deals with either the non-characteristic boundary case or the boundary characteristic case where the characteristic eigenvalue vanishes identically. The boundary condition is assigned by imposing that 
\be \label{e:bcb}
   \mathbf b (\mf v(t, 0) ) = \mf h (t) \quad \text{a.e. $t \in \R_+$}.
\eq
In the previous expression, when the domain is $x>0$ the boundary datum $\mf h$ attains values in $\mathbb R^{p}$, where $p$ is the number of strictly positive eigenvalues of the jacobian matrix of $\mf D\mf f$. The function $\mf b: \R^N \to \R^p$ is given and must sastify suitable assumptions. Once again, it is worth stressing that the assumption that the boundary characteristic eigenvalue vanishes identically is foundamental here: in the general boundary characteristic case where one eigenvalue can attain both positive and negative eiegenvalues, the number of boundary conditions to prescribe at the boundary is not determined \emph{a priori} as it depends on the sign of the characteristic eigenvalue of the solution at the boundary, and hence one cannot formulate the boundary condition as in~\eqref{e:bcb}. 
It should be also noted that all the above works~\cite{Goodman,LiuIBVP,NishidaSmoller,SableT} require suitable smallness assumptions on the total variation of the data (or on some related quantities, as for instance in~\cite{NishidaSmoller}) and in several cases the analysis extends to equations defined on the strip $0 < x < 1$. In particular, in~\cite{SableT} the author introduces an \emph{ad-hoc} condition to prevent the possible amplification of reflected waves in the case.  
\subsection{Results \emph{via} wave front-tracking algorithms}
In~\cite{Amadori} Amadori established global existence results for admissible solution of~\eqref{e:claw} defined on the domain $x>0$. The proof requires, as usual, smallness assumptions on the data and relies on the introduction of a suitable wave front-tracking algorithm. In~\cite{Amadori} the author considers both non characteristic and 
characteristic boundaries, and in the latter case there is no constraint on the characteristic eigenvalue, which can attain both positive and negative values. In the non characteristic boundary case, the boundary condition is prescribed as in~\eqref{e:bcb}, whereas in the boundary characteristic case Amadori relies on a formulation introduced in a previous work by Dubois and LeFloch~\cite{DuboisLeFloch} which we now briefly discuss. Given
$\mf v_b: \R_+ \to \R^N$, in~\cite{Amadori} one assigns the boundary condition in the boundary characteristic case 
by prescribing that 
$$
    \mf v(t, 0) \sim_{\ast} \mf v_b (t) \quad \text{a.e. on $\R_+$},
$$
where the relation $\sim_{\ast}$ is defined as follows. Consider the Riemann problem obtained by coupling~\eqref{e:claw} with the initial datum
$$
   \mf v(0, x) =
  \left\{ 
  \begin{array}{ll}
   \mf u^- & x<0 \\
   \mf u^+ & x>0 \\ 
  \end{array}
  \right.
$$ 
and set $\mf u^-: = \mf v_b (t)$,  $\mf u^+:=\mf v (t, 0)$ (the trace of the solution $\mf v$). Then $ \mf v(t, 0) \sim_{\ast} \mf v_b (t)$ if the solution of the Riemann problem only contains waves that are not entering the domain, namely have non-positive speed. It is worth stressing that this way of assigning the boundary condition does not involve any information on the underlying viscous mechanism and in general it is different from the one provided by Definition~\ref{d:equiv}. To appreciate this difference, and henceforth the difference between the relations $\sim_\ast$ and $\sim_{\mf D}$, it suffices to consider the linear case~\eqref{e:linear}. The boundary condition considered in~\cite{Amadori,DuboisLeFloch} dictates that, for a.e. $t \in \R_+$, $\mf v_b(t) - \mf v(t, 0)$ belongs to the generalized eigenspace of the matrix $\mf A$ associated to non-positive eigenvalues, whereas prescribing~\eqref{e:ibvp} we prescribe that, for a.e. $t \in \R_+$, $\mf v_b(t) - \mf v(t, 0)$ belongs to the generalized eigenspace of the matrix $\mf D^{-1}\mf A$ associated to non-positive eigenvalues. Since in general the eigenspaces of $\mf A$ and $\mf D^{-1}\mf A$ do not coincide, the two conditions are different. 
\subsection{Vanishing viscosity approximation}
To conclude this brief overview, we mention that, in the case of the \emph{artificial viscosity} where $\mf D$ is the identity matrix, there are some results establishing convergence of the vanishing viscosity approximation~\eqref{e:vclaw} in domains with boundaries, see~\cite{AnconaBianchini,ChenFrid,Spinolo}. As a byproduct, these results yield global existence for solutions of the conservation law~\eqref{e:claw}.  We stress once more, however, that in the case of initial-boundary value problems  global-in-time convergence proofs are limited to the artificial viscosity case, and hence in general the limit does not provide the physically relevant solution. Convergence results for the physical viscosity have been only obtained for regular data and on small time intervals (typically, before the classical solution of the conservation law breaks down), see for instance the works by Gisclon~\cite{Gisclon} and by Rousset~\cite{Rousset}. See also Joseph and LeFloch~\cite{JosephLeFloch} for the convergence of a  different self-similar viscous approximation in the case of Riemann-type data.   

\section{Some technical comments on the proof of Theorem~\ref{t:main}} 
\label{s:technical}
In this paragraph we provide some handwaving comment on the main technical novelties of the construction in~\cite{AFL}. In \S\ref{ss:gc} we make some general comment on the proof of Theorem~\ref{t:main} and mention that the most innovative technical points, which we discuss in \S\ref{ss:ie}, are a very detailed boundary analysis, and some new interaction estimates. In \S\ref{ss:wft} to frame the interaction estimates in the right context  we briefly and informally overview the wave front-tracking construction for the Cauchy problem. 

\subsection{General comments} \label{ss:gc} 
As mentioned in \S\ref{ss:over}, the proof of Theorem~\ref{t:main} relies on the introduction of a new wave front-tracking algorithm. At a very high level, the main novelty of the scheme discussed in~\cite{AFL} is that it incorporates detailed information on the underlying viscous approximation and in particular on the transient behavior  as $\ee \to 0^+$ of~\eqref{e:vclaw} at the domain boundary: for instance, we introduce a new Glimm-type interaction functional that involves, among other things, the strength of the center component of the boundary layer sitting at the domain boundary. Note that incorporating information on the viscous mechanism is possible owing to the detailed analysis of the structure of the boundary layer profiles (that is, of the solutions of~\eqref{e:bl}) that was done in~\cite{BianchiniSpinoloARMA,BianchiniSpinolo}. 

At a lower and more detailed level, it is worth highlighting two main technical novelties introduced in~\cite{AFL}. The first one is a new interaction estimate that applies when the boundary characteristic field is genuinely nonlinear and a wave of the characteristic family hits the boundary. We comment on this estimate in the following, for the time being we point out that the boundary characteristic field is genuinely nonlinear if, for instance, system~\eqref{e:claw} are the compressible Euler equations written in Eulerian coordinates and the modulus of the fluid velocity is close to the sound speed.  The second technical novelty that we want to single out is a detailed analysis of the behavior of the wave front-tracking approximation at the domain boundary, see~\cite[\S9]{AFL}. This analysis is instrumental in showing that the limit of our wave front-tracking approximation satisfies indeed the boundary condition~\eqref{e:ibvp}, and exploits fine properties of the wave front-tracking algorithm that are discussed in~\cite[Ch.10]{Bressan}. 
\subsection{The wave front-tracking algorithm} \label{ss:wft}
 Prior to discussing  the main technical novelties that we touched upon before, we have to briefly recall some very basic steps in the construction of the wave front-tracking algorithm, and we refer to the fundamental references~\cite{Bressan,Dafermos,HoldenRisebro} for a complete discussion. In a nutshell, any wave front-tracking scheme provides an approximation of the admissible solution of~\eqref{e:claw} through piecewise constant functions with a finite number of discontinuity lines, the so-called wave fronts after which the algorithm is named. In the case of the Cauchy problem, the first step of the algorithm is constructing a piecewise constant approximation of the initial datum. At every discontinuity point of the approximate initial datum, one constructs a piecewise constant solution of the Riemann problem with data given by the left and right limits of the approximate initial datum at the discontinuity point. In particular, to obtain a piecewise constant solution one must approximate rarefaction waves through a wave fan involving finitely many states separared by discontinuity lines. Once the approximation of the solution of the Riemann problem at each single discontinuity is determined, one can construct a local-in-time approximate solution of the Cauchy problem by juxtaposing the approximate solutions of the Riemann problems.  This solution can be extended in time until two distinct wave fronts cross each other: to further extend the solution, one approximately solves the Riemann problem determined at the interaction point, and in this way defines an approximate solution that is defined up to the second interaction time, where one approximately solves the corresponding Riemann problem, and so on. 
Needless to dwell on the details, two of the main challenges in constructing global-in-time approximate solutions with a wave front-tracking algorithm are on the one hand to prevent the formation of infinitely many wave fronts in finite time, and on the other to control the total variation growth. To tackle the first challenge, a by-now-standard technique, which we also use in~\cite{AFL}, is the introduction of so-called \emph{non-physical} fronts: very loosely speaking, wave fronts that are somehow negligible are forced to travel with very high speed and to have minimal interactions with the ``physical" fronts, see~\cite[Ch. 7]{Bressan} for the precise construction. To tackle the second issue and find a uniform bound on the total variation of the approximate solution a key point is establishing the so-called interaction estimates, namely extracting detailed  information on the strength of the waves generated when two wave fronts interact, in terms of the strength of the interacting waves. As a matter of fact, establishing suitable interaction estimates is one of the key points of the constrution of the wave front-tracking algorithm, and the precise form of these estimate determines the
specific form of the Glimm-type functionals used to control the total variation. 
\subsection{Main technical novelties of~\cite{AFL}}\label{ss:ie}
In the case of the initial-boundary value problem~\eqref{e:claw},\eqref{e:ibvp}, the first step in our construction is establishing piecewise constant approximations of the initial datum $\mf v_0$ and of the boundary datum $\mf v_b$. At every discontinuity point in the approximate boundary datum we approximately solve the so-called boundary Riemann problem, see~\cite[\S4]{AFL} for the precise construction and in particular \S4.1 for an heuristic overview of the basic ideas underpinning the analysis. Without dwelling on the details, what is worth pointing out here is that our approximate solution of the boundary Riemann problem relies on the exact solution discussed in~\cite{BianchiniSpinoloARMA,BianchiniSpinolo}, and provides precise information on the boundary layer profiles satisfying~\eqref{e:bl}, and in particular on the size of their component on a suitable center manifold. In the following, we term this size $\xi_k$, with $k$ denoting the index of the boundary characteristic family\footnote{in other words, the $k$-th eigenvalue of $\mf D \mf f (\mf u)$ is the one that can vanish}, and as we will see it will appear in our main interaction estimate.

\begin{figure} \label{f:1}
\begin{center}
\begin{tikzpicture}
\draw[line width=0.4mm,->]   (0, -1) -- (0, 6) node[anchor=west] {$t$};
\draw[line width=0.4mm,red] (1, -1) -- (0, 2);
\draw (0, 2) node[anchor=east] {$\tau$};
\draw[line width=0.4mm,red] (0, 2) -- (1, 6);
\draw[line width=0.4mm,blue] (0, 2) -- (4, 6);
\draw[line width=0.4mm,blue] (0, 2) -- (5, 5.5);
\draw[line width=0.4mm,magenta] (0, 2) -- (9, 5);
\end{tikzpicture}
\caption{A wave front of the boundary characteristic family (red) hitting the boundary at time $\tau$, and the waves entering the domain that are generated at the interaction.}
\end{center}
\end{figure}
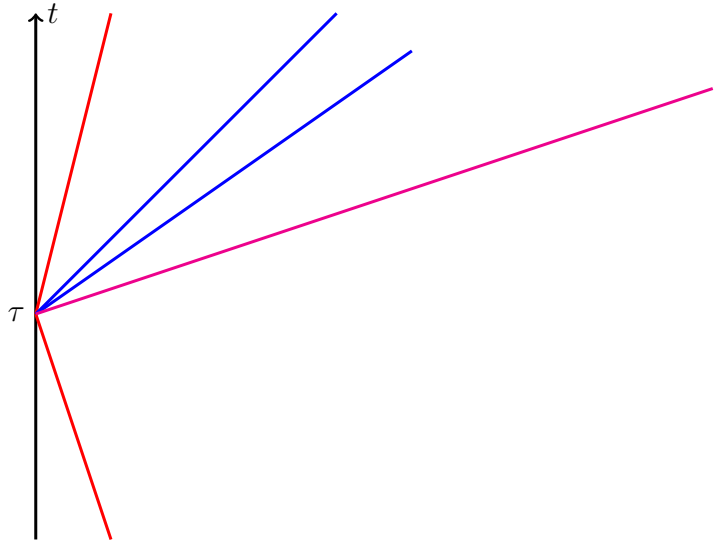

As in the case of the Cauchy problem, for the initial-boundary value problem by juxtaposing approximate solutions of Riemann and boundary Riemann problems we can define a local-in-time solution which can be extended up to the first interaction time. Note that, in the case of the initial-boundary value problem, interactions occur at times where either two wave fronts collide, or a wave front hits the boundary. To handle the first case we basically proceed as in the case of the Cauchy problem and employ the interactions estimates established in~\cite[Ch 7]{Bressan}. The non standard part is in the handling of the second type of collision. In what follows we focus on the most innovative part of our contribution and assume that the boundary characteristic field is genuinely nonlinear. We consider a wave front of the boundary characteristic family that hits the domain boundary $x=0$ at the time $\tau$, see Figure~\ref{f:1} for a representation. 

Our main interaction estimate is equation (6.5) in~\cite{AFL} and states that there is a constant $C$, only depending on system~\eqref{e:vclaw} and on the value $\mf v^\ast$ in the statement of Theorem~\ref{t:main}, such that the following holds: if we term $|\Delta V(\tau)|$ the increase of the total variation of the approximate solution at time $\tau$  we have 
\be \label{e:mainie}
   |\Delta V(\tau)| \leq C |s_k| ([\varsigma_k]^- + |\xi_k|).
\eq
In the previous expression, $s_k$ denotes the size of the hitting wave, and $\xi_k$ is either the size of the center component of the boundary layer sitting at the domain boundary before the interaction, or the size of the $0$-speed Lax admissible shock located at $x=0$ before the interaction. As a matter of fact, the cases where the boundary layer has a center component and where there is a Lax shock tangent at the boundary are mutually exclusive owing to the genuine nonlinearity assumption, see~\cite[\S4.3]{AFL} for the details of the analysis. Finally, $[\varsigma_k]^-$ denotes the negative part of the quantity $\varsigma_k$, which in turn is the shock speed if the hitting front is a shock wave, and the eigenvalue of the boundary characteristic family evaluated at the leftmost state if the hitting front is a rarefaction wave. To the best of our knowledge,~\eqref{e:mainie} is the first boundary interaction estimate involving the speed of the hitting wave, as well as the first one involving the size of the boundary layer sitting at the domain boundary. From the technical standpoint, the proof of~\eqref{e:mainie} is rather delicate and it is complicated by the fact that the functions 
involved in the definition of the solution of the boundary Riemann problem (and in particular the so-called \emph{characteristic wave fan curve of admissible states} defined in~\cite[\S4.3]{AFL}) are only Lipschitz continuous and fail to attain $C^1$ regularity. 

After establishing our interaction estimates, in~\cite{AFL} we introduce a new Glimm-type functional, which in turn allows us to control the total variation growth and hence establish strong $L^1_{\mathrm{loc}}$ compactness of the wave front-tracking approximation \emph{via} the Helly-Kolmogorov Compactness Theorem. 
What is then left to prove is that any accumulation point of the wave-front tracking approximation is a solution of the initial-boundary value problem~\eqref{e:claw},\eqref{e:ibvp}. If some parts of the limit analysis are by now fairly standard and rely on the well-established techniques discussed in~\cite{Bressan}, showing that the any accumulation point satisfies the boundary condition in the sense of Definition~\ref{d:equiv} is a highly nontrivial point whose proof requires an accurate analysis of the boundary behavior of the approximate solutions and relies on the powerful techniques discussed in~\cite[Ch. 10]{Bressan}. 
We now touch upon some technical details, and to highlight the heart of the matter  we focus on the case of a genuinely nonlinear boundary characteristic field, which has a richer and more interesting behavior.  Heuristically speaking, there are two basic mechanisms that could in principle prevent the limit of a sequence of wave front-tracking approximation from attaining the boundary condition in the sense of Definition~\ref{d:equiv}: on the one hand, in the approximation there might be a sequence of shocks with vanishing speed that accumulate at the domain boundary. On the other, in the approximation there might a sequence of boundary layers with nontrivial center components, that in the limit converge to a $0$-speed Lax shock. In both cases, and at the price of a fairly technical analysis, in~\cite{AFL} we show that the limit satisfies  Definition~\ref{d:equiv} with a non-trivial $0$-speed shock sitting at the domain boundary. This is very looosely speaking one of the two main ingredients of the boundary analysis in~\cite{AFL}, the other one being a compactness result for the boundary trace of the flux functions of the wave front-tracking approximations, see~\cite[Lemma 9.2]{AFL}. This is a nontrivial result since we expect that the trace of the solutions may instead have an higly oscillatory behavior, see for instance the related  example in~\cite[\S4.3]{DMS}.

\section*{Acknowledgments}
LVS wishes to thank Professors Mikhail Feldman, Dehua Wang, Wei Xiang and Tong Yang for the invitation to contribute to the celebration of Professor Gui-Qiang Chen's birthday. All authors are members of the GNAMPA group of INDAM and of the PRIN Project 20204NT8W4 (PI Stefano Bianchini). FA and LVS are also members of the PRIN 2022 PNRR Project P2022XJ9SX (PI Roberta Bianchini) and LVS is a member of the PRIN Project 2022YXWSLR (PI Paolo Antonelli). 
\bibliographystyle{plain}
\bibliography{wft}
\end{document}